\documentstyle{article}

\begin{document}

\hfill{q-alg/9610010}

\hfill{October, 1996}

\vspace{10mm}

\centerline{\bf THE QUANTUM SYMPLECTIC CAYLEY-KLEIN GROUPS}
\vspace*{0.37truein}
\centerline{\footnotesize N.A. GROMOV, I.V. KOSTYAKOV and V.V. KURATOV}
\vspace*{0.015truein}
\centerline{\footnotesize\it Department of Mathematics, Komi Science Centre,}
\baselineskip=10pt
\centerline{\footnotesize\it Ural Division, Russian Academy of Sciences,}
\baselineskip=10pt
\centerline{\footnotesize\it Syktyvkar, 167000, Russia }
\baselineskip=10pt
\centerline{\footnotesize\it E-mail : parma@omkomi.intec.ru}
\vspace*{0.21truein}

\abstract {The contraction method
applied to the construction of the nonsemisimple quantum
symplectic Cayley-Klein groups $ Fun(Sp_q(n;j)) $.
 This groups has been realised
as Hopf algebra of the noncommutative functions over the algebra
with nilpotent generators.  The  dual quantum algebras $ sp_q(n;j) $
are constructed.}

\textheight=7.8truein
\textwidth=5truein
\renewcommand{\thefootnote}{\alph{footnote}}

\section{Introduction}
\noindent
A description of nonsemisimple quantum symplectic
groups and algebras is suggested.
It is known, that the groups  motion of Cayley-Klein space,
can be  described uniformly \cite{1}, if  the parameter
$  j ,(j=1,\iota ,$ where $ \iota\not=0 ,{\iota}^{2}=0 ) $
is introduced to the group.
The symplectic groups is obtained by this method will be denoted as
$ Sp(n;j), $
and this introduction of parameter $ j $ shall be called
the Cayley-Klein structure introduction  on the group.
The Cayley-Klein structure on the quantum symplectic group is defined and
the description of nonsemisimple quantum symplectic groups is obtained.
We used the notation \cite{2}.

\section{ Symplectic algebra sp(n;j) }
\noindent
Consider the space $ {\rm R^n}(j) \otimes {\rm R^n}(j) $, which is obtained
from $ 2n$-dimensional Euclidean space $ {\rm R^n} \otimes {\rm R^n} $
by mapping
$$
 \psi: \enskip {\rm  R^n \otimes R^n \rightarrow R^n}(j) \otimes
 {\rm R^n }(j) ,
$$
$$
\psi  x_k = x_kJ_k, \quad  \psi x_{k^{\prime}} = x_{k^{\prime}}J_k, \quad
k^{\prime} = 2n + 1 - k, \enskip
J_k= \prod_{m=1}^{k-1} j_m ,
\enskip k=1,2 \ldots n .
$$
Group $ Sp(n;j) $ is defined as the group of transformation of
 $ 2n-$dimensional space
$ R^n(j) \otimes R^n(j) $,
preserving the bilinear form
\begin{equation}
(x,y)= \sum_{k=1}^{n}(x_ky_{k^{\prime}}-x_{k^{\prime}}y_k)J_k^2.
\end{equation}
Cartesian  coordinates $ x_k, y_k, \enskip k=1,2 \ldots n, $
belongs to the first, and  $ x_{k^{\prime}},y_{k^{\prime}} $---
to the second factor in the direct product of spaces.
For $ j $-structure to be preserved, parameter $ j $ must be introduced
in matrix of transformation $ T(j) = (T_{ij})^n_{ i,j = 1} $ as follows
$$
 T_{ij} =  J^k_m t_{ij} ,\  where \enskip  k = i \enskip if \enskip i \leq n;
 \quad k = i^{\prime}, \enskip i > n;
 \quad m = j,  \enskip j \leq n;
$$
$$
 m = j^{\prime} , \enskip j > n, \qquad
 J^k_m = \prod_{n=min(k,m)}^{max(k,m)-1} j_n .
$$
The matrix elements of $ T $ satisfies the additional relation
  \begin{equation}
     T^tC_0T=C_0.
  \end{equation}
The elements of matrix $ C_0 $ are as follows
$ (C_0)_{ij} = \varepsilon_i \delta_{ij} , \enskip
 \varepsilon_i = 1 $ if $ i=1,2 \ldots n $ and
$ \varepsilon_{i} = -1 $ for $ i = n-1, \ldots , 2n $.
The general element of the Shevalley basis of symplectic algebra
$ sp(n;j) $ appears as
\begin{equation}
 H= \left (\begin{array}{cccccc}
   h_{1}    & j_1 x_1^+ &  \cdot      & \cdots    & \cdots  & \cdot      \cr
j_1 x_1^-   & h_2       &  j_2 x_2^+  & \cdots    & \cdots  & \cdot      \cr
   \cdot    & j_1 x_2^+ & \cdots      & \cdots    & \cdots  & \cdot      \cr
   \cdot    & \cdots    & \cdots      & \cdots    & \cdots  & \cdots     \cr
   \cdot    & \cdots    & \cdots      &-j_2 x_2^- & - h_1   & -j_1 x_1^+  \cr
   \cdot    & \cdots    & \cdots      & \cdots    & - j_1 x_1^- & - h_2
            \end{array} \right ) .
\end{equation}
Choosing the generators of algebra $ sp(n;j) $ in the following form
 \begin{eqnarray}
 H_i &  = &  e_{ii} - e_{ i^{\prime}i^{\prime}} - e_{i+1,i+1}
 +  e_{(i+1)^{\prime},(i+1)^{\prime}} ,\quad
 H_n = e_{nn} - e_{n^{\prime},n^{\prime}} ,    \nonumber \\
 X_i^{+} &  = & j_i(e_{i,i+1} - e_{(i+1)^{\prime},i^{\prime}})  ,\quad
 X_i^{-} = j_i(e_{i+1,i} - e_{(i)^{\prime},(i+1)^{\prime}} )  ,   \nonumber \\
 X_n^{+} & = & j_n  e_{n,n+1} ,\quad
 X_n^{-} = j_n e_{n+1,n} ,  \quad i = 1, \ldots , n-1 ,
 \end{eqnarray}
where $ (e_{ij})_{km}=\delta_{ik}\delta_{jm}, $
the cÏmmutation relation for the Shevalley basis $ sp(n;j) $ looks as
 \begin{equation}
 \left [H_i,X_j^{ \pm} \right ]   =  \pm A_{ij} X_j^{ \pm} ,\quad
 \left [X_i^{+},X_j^{-} \right ]  =  \delta_{ij} j_i^2 H_j , \quad
 \end{equation}
where $ A_{ij} $ is a Cartan matrix ---
$ A_{ii} = 2, \enskip  A_{i,i-1} = A_{i-2,i-1} = -1, \enskip
A_{n-1,n} = - 2 $.

\section{Cayley-Klein structure on quantum symplectic groups
 $ Fun(Sp_q(n;j)) $ }
\noindent
The Cayley-Klein structure is introduced to
quantum symplectic groups \\ $ Fun(Sp_q(n)) $ ,
as well as in the classical case with the help of the map
\begin{eqnarray}
 \psi \hat x_k = x_kJ_k, &   \psi \hat x_{k^{\prime}} = x_{k^{\prime}}J_k,&
 k = 1, \ldots , n,
\end{eqnarray}
where $ \hat{x}_i $ forms the basis of quantum N-dimension
 symplectic Cayley-Klein space
    $ Sp_q^N( \rm C;j) $ with commutation relations
 \begin{equation}
   \hat R_q( \hat x \otimes \hat x)=q(\hat x \otimes \hat x).
 \end{equation}
Parameter $ j $  is introduced in matrix $ T $ similarly to the classical case.
Commutation and additional relations for generators $ t_{ik} $ are defined by
  \begin{equation}
  RT_1(j)T_2(j)=T_2(j)T_1(j)R ,
  \end{equation}
  and
  \begin{equation}
     T(j)CT^t(j)C^{-1}=CT^t(j)C^{-1}T(j)=I ,
  \end{equation}
whith \ $ T_1(j) = T(j) \otimes I ,
\quad T_2(j)=I \otimes T(j) , \quad
 C = C_0 q^{\rho} , \enskip \\
\rho = {\rm diag}(n,n-1, \ldots,1,-1, \dots ,-n) $.
The coproduct, counit and antipode is defined generally as
$$
 \triangle(T(j))=T(j)\stackrel{.}{\otimes}T(j), \quad
 \epsilon(T(j))=I,
$$
  \begin{equation}
  S(T(j))= CT^t(j)C^{-1} .
  \end{equation}
  The  "action"
  $ \delta^{\prime} $:
  $ Sp_q^N( {\rm C;j}) \rightarrow Fun(Sp_q(n;j)) \otimes
  Sp_q^N({ \rm C;j})
  $ of quantum group \\
  $ Fun(Sp_q(n;j)) $ on quantum symplectic Cayley-Klein space
  $ Sp_q^N({ \rm C };j) $,
 is defined by
  $ \delta^{\prime} = T(j) \stackrel{.}{ \otimes} \hat x $,
 preserve the bilinear form
   $\hat x^{t} \stackrel{.}{\otimes} C \hat x $.

\section{Quantum algebra $ sp_q(n;j) $ as dual to $ Fun(Sp_q(n;j)) $ }
\noindent
 Dual to $ Fun(Sp_q(n;j)) $
 the quantum algebra $ sp_q(n;j) $ is
  defined by the relation \cite{2}
$$
   < L^{(\pm)}(j),T(j)>= R^{(\pm)}.
$$
 \begin{equation}
  R^{(-)}  =   R_{q}^{-1} = R_{q^{-1}}, \quad
  R^{(+)}  =   PR_{q}P ,\quad Pu \otimes v = v \otimes u.
  \end{equation}
 Generators of the algebra $ sp_q(n;j) $
 have the following form
$$
\left (  L^{(+)}(j) \right ) _{ij} = \tilde J^k_m l_{ij}^{(+)} ,\quad for
 \enskip i \leq j \enskip and
\enskip \left (  L^{(+)}(j) \right ) _{ij} =  0  \quad for \enskip i > j ,
$$
$$
\left (  L^{(-)}(j) \right) _{ij} = \tilde J^k_m l_{ij}^{(-)} ,\quad  for
 \enskip
i \geq j \enskip and \enskip
\left (  L^{(-)}(j) \right) _{ij} =   0  \quad for \enskip i < j,
$$
$
 \qquad k = i , \enskip i \leq n;
 \quad k = i^{\prime}, \enskip i > n;
 \quad m = j , \enskip j \leq n;
 \quad m = j^{\prime} , \enskip j > n,
$
$$
\tilde J^k_m = \prod_{n=min(k,m)}^{max(k,m)-1} j_n^{-1}.
$$
Commutation and additional relations for generators of    $ sp_q(n;j) $
algebras are given by
 \begin{eqnarray}
   R^{(+)}L_1^{(\varepsilon)}(j)L_2^{(\varepsilon)}(j) & =
    & L_2^{(\varepsilon)}(j)L_1^{(\varepsilon)}(j)R^{(+)} ,   \\
   R^{(+)}L_1^{(+)}(j)L_2^{(-)}(j) & = & L_2^{(-)}(j)L_1^{(+)}(j)R^{(+)}, \\
    L(j)C^tL^{t}(j)(C^{-1})^{t}         &         =         &
 C^{t}L^{t}(j)(C^{-1})^{t}L(j)=I , \label{4}
 \end{eqnarray}
where
 $
  L_1^{( \varepsilon )}(j)=L^{( \varepsilon )}(j) \otimes I, \quad
  L_2^{( \varepsilon)}(j)=I \otimes L^{( \varepsilon )}(j), \quad
  \varepsilon=\pm
 $.
 Coproduct and counit are defined  generally as
  \begin{equation}
  \triangle(L^{\pm}(j))=L^{\pm}(j)\stackrel{.}{\otimes}L^{\pm}(j),\quad
  \epsilon(L^{\pm}(j))=I.
  \end{equation}

The isomorphism between algebra   $   sp_{q}(n;j) $
and universal enveloping algebra  $  U_z (sp( n;j)) $
is defined by the following expression:
\begin{equation}
   l_{ii}^{(+)}=e^{z \tilde H_i} ,\quad   l_{nn}^{(+)}=e^{z \tilde H_n },
\end{equation}
\begin{equation}
l_{i, i+1}^{(+)}=\lambda q^{-\frac{1}{2} } X^+_i e^{ \frac{z}{2}(\tilde H_i +
\tilde H_{i+1}) },
\quad   l_{i+1, i}^{(-)}= - \lambda q^{\frac{1}{2} } X^{-}_{i}
 e^{- \frac{z}{2} (\stackrel{ \sim} {H_i} +
 \tilde H_{i+1}) }
 \end{equation}
   \begin{equation}
   l_{n,n+1}^{(+)}=\lambda  X^+_n e^{\frac{z \tilde H_n}{2}},\quad
   l_{n+1,n}^{(-)}=\lambda  X^-_n e^{-\frac{z \tilde H_n}{2}},
  \end{equation}
  $$
  H_i= \tilde H_i - \tilde H_{i+1}, \quad
  H_n = 2 \tilde H_n .
  $$
Commutation relations in $  U_z(sp(n;j)) $ are given by
 \begin{eqnarray}
 [H_i,H_j ]=0, &   [H_i,X_j^{  \pm }]= \pm ( \alpha_i,
  \alpha_j ) X_j^{ \pm }, &  \nonumber \\
 \left [X_i^+,X_j^- \right]= j^{2}_i [H_i]_q \delta_{ij} & [H_i]_q =
 \frac {sh(zH_i)}{sh(z)}. &  \label{1}
 \end{eqnarray}
 Equations (~\ref{1})  defines the structure of universal enveloping algebra
  $  U_z(sp(n;j)) $.
   In the case $ j_i = \iota_i  $
   we have obtained a new nonsemisimple algebra.

\section{ Case $ n = 2 $ as example }
\noindent
The quantum symplectic group $ Fun(Sp_q(2;j)) $
is generated by the matrix with noncommutative generators
\begin{equation}
 T(j)= \left (
 \begin{array}{cccc}
  t_{11} & jt_{12} & jt_{13} &  t_{14} \cr
 jt_{21} &  t_{22} &  t_{23} & jt_{24} \cr
 jt_{31} &  t_{32} &  t_{33} & jt_{34} \cr
  t_{41} & jt_{42} & jt_{43} &  t_{44}
 \end{array} \right ) .
\end{equation}
For example in explicit form the antipode is as follows
\begin{equation}
 S(T(j))= \left (\begin{array}{cccc}
             t_{44}     & jq^{-1}t_{34} & -jq^{-3}t_{24} & -q^{-4}t_{14} \cr
             jqt_{43}   & t_{33}        & -q^{-2}t_{23}  & -jq^{-3}t_{13} \cr
          -jq^3t_{42}   & -q^2t_{32}    & t_{22}         & jq^{-1}t_{12} \cr
           -q^4t_{41}   & -jq^3t_{31}   & jqt_{21}       & t_{11}
                 \end{array} \right ).
\end{equation}
The generators of the dual to $ Fun(Sp_q(2;j)) $ quantum algebra $ sp_q(2;j) $
written in compact matrix form as
  \begin{equation}
             L^{(+)}(j) = \left (\begin{array}{cccc}
         l_{11}^{(+)} & j^{-1}l_{12}^{(+)} & j^{-1}l_{13}^{(+)} & l_{14}^{(+)}   \cr
             0        & l_{22}^{(+)}  & l_{23}^{(+)}  & j^{-1}l_{24}^{(+)}  \cr
             0        &      0    & l_{33}^{(+)}  & j^{-1}l_{34}^{(+)}  \cr
             0        &      0    &     0     & l_{44}^{(+)}
              \end{array} \right ) ,
  \end{equation}
  \begin{equation}
 L^{(-)}(j) = \left (\begin{array}{cccc}
              l_{11}^{(-)}      &    0      &     0     &    0      \cr
            j^{-1}l_{21}^{(-)}  & l_{22}^{(-)}  &     0     &    0      \cr
            j^{-1}l_{31}^{(-)}  &  l_{32}^{(-)} &  l_{33}^{(-)} &    0      \cr
    l_{41}^{(-)} & j^{-1}l_{42}^{(-)} & j^{-1}l_{43}^{(-)}  &
 l_{44}^{(-)}
              \end{array} \right )  .
   \end{equation}

Commutation and additional relations for generators of $ sp_q(2;j) $
algebras are given by
$$
 ql_{11}^{(+)}l_{12}^{(+)} = l_{12}^{(+)}l_{11}^{(+)},\quad
  l_{11}^{(+)}l_{23}^{(+)} = l_{23}^{(+)}l_{11}^{(+)},\quad
  l_{11}^{(+)}l_{32}^{(-)} = l_{32}^{(-)}l_{11}^{(+)},
$$
$$
  l_{11}^{(+)}l_{21}^{(-)} = ql_{21}^{(-)}l_{11}^{(+)},\quad
  l_{11}^{(\pm)}l_{kk}^{(\pm)} = l_{kk}^{(\pm)}l_{11}^{(\pm)},\quad k=1,
 \dots ,4 ,
$$
$$
  ql_{22}^{(+)}l_{21}^{(-)} = l_{21}^{(-)}l_{22}^{(+)},\quad
  l_{22}^{(+)}l_{32}^{(-)} = q^2l_{32}^{(-)}l_{22}^{(+)},\quad
  l_{22}^{(+)}l_{12}^{(+)} = ql_{12}^{(+)}l_{22}^{(+)},
$$
$$
  q^2l_{22}^{(+)}l_{23}^{(+)} = l_{23}^{(+)}l_{22}^{(+)},\quad
 ql_{23}^{(+)}l_{12}^{(+)} - l_{12}^{(+)}l_{23}^{(+)}=\tilde \lambda
 l_{13}^{(+)}l_{22}^{(+)},
\quad \tilde \lambda =q^2-q^{-2},
$$
$$
  [l_{12}^{(+)},l_{21}^{(-)}]=\lambda(l_{11}^{(+)}l_{22}^{(-)}-
 l_{11}^{(-)}l_{22}^{(+)}),
  \quad
 [l_{23}^{(+)},l_{32}^{(-)}]=\tilde \lambda (l_{22}^{(+)}l_{22}^{(+)}-
 l_{22}^{(-)}l_{22}^{(-)}),
$$
\begin{equation}
 ql_{32}^{(-)}l_{21}^{(-)} - l_{21}^{(-)}l_{32}^{(-)} = \tilde \lambda
 l_{31}^{(-)}l_{22}^{(-)},\quad
  l_{12}^{(+)}l_{32}^{(-)} = ql_{32}^{(-)}l_{12}^{(+)},\quad
  ql_{23}^{(+)}l_{21}^{(-)} = l_{21}^{(-)}l_{23}^{(+)} .
\end{equation}
Antipode :
  \begin{equation}
           S(L^{(+)}(j)) = \left (\begin{array}{cccc}
   l_{44}^{(+)} & j^{-1}ql_{34}^{(+)} & -j^{-1}q^{3}l_{24}^{(+)}
 & -q^4l_{14}^{(+)}        \cr
   0      &     l_{33}^{(+)}     &     -q^2l_{23}^{(+)}     &
 -j^{-1}q^3l_{13}^{(+)}  \cr
     0  &   0   & l_{22}^{(+)}        & j^{-1}ql_{12}^{(+)}     \cr
     0  &  0        &     0               & l_{11}^{(+)}
              \end{array} \right ) ,
  \end{equation}
  \begin{equation}
       S(L^{(-)}(j)) = \left (\begin{array}{cccc}
    l_{44}^{(-)} & 0 & 0 & 0 \cr  j^{-1}q^{-1}l_{43}^{(-)}  &
    l_{33}^{(-)} & 0 & 0 \cr
   -j^{-1}q^{-3}l_{42}^{(-)}  & -q^{-2}l_{32}^{(-)} &
 l_{22}^{(-)} &    0      \cr
   -q^{-4}l_{41}^{(-)}
 &  -j^{-1}q^{-3}l_{31}^{(-)}  &  j^{-1}q^{-1}l_{21}^{(-)}  &
 l_{11}^{(-)}
              \end{array} \right )  .
   \end{equation}

The isomorphism between quantum algebra $ sp_{q}(2;j) $
and universal enveloping algebra $ U_z (sp(2;j)) $
is defined by the following expressions
 ( $ q=e^{ \varepsilon z},\enskip \varepsilon= \pm1 $ ):
\begin{equation}
   l_{11}^{(+)} = e^{z \tilde H_i} ,\quad l_{22}^{(+)} = e^{z \tilde H_n },
\end{equation}
\begin{equation}
l_{12}^{(+)}=\lambda q^{-\frac{1}{2} } X^+_1 e^{ \frac{z}{2}(\tilde H_1 +
\tilde H_2)},
\quad l_{21}^{(-)} = - \lambda q^{\frac{1}{2} } X^{-}_{1} e^{- \frac{z}{2}
 (\stackrel{ \sim} {H_i} +
\tilde H_2) },
\end{equation}
   \begin{equation}
   l_{23}^{(+)}=\lambda q^{-\frac{1}{2} } X^+_2 e^{
 \frac{z \tilde H_2}{2}},\quad
 l_{32}^{(-)}=\lambda  q^{\frac{1}{2}  }  X^-_2   e^{-\frac{z
 \tilde H_2}{2}},
  \end{equation}
  $$
  H_1= \tilde H_1 - \tilde H_2, \quad
  H_2 = 2 \tilde H_2 .
  $$
Other generators may be obtained from additional relation (~\ref{4})
 \begin{eqnarray*}
 l_{13}^{(+)} & = &
 ( \tilde \lambda )^{-1}(q l_{23}^+l_{12}^{(+)} -
 l_{12}^{(+)}l_{23}^{(+)})l_{22}^{(-)}   , \quad
  l_{31}^{(-)} =  (\tilde \lambda )^{-1}
  (ql_{32}^{(-)}l_{21}^{(-)}-l_{21}^{(-)}l_{32}^{(-)} ) l_{22}^{(+)}, \\
l_{14}^{(+)} & = & j^{-2} \lambda^{-1}(l_{13}^{(+)}l_{12}^{(+)}
 - q^{2}l_{12}^{(+)}l_{13}^{(+)})l_{11}^{(-)} ,\qquad
  l_{34}^+  =  -l_{12}^{(+)}l_{11}^{(-)}l_{22}^{(-)}, \\
  l_{41}^{(-)} & = & j^{-2} \lambda^{-1}(l_{31}^{(-)}l_{21}^{(-)}
 - q^{2}l_{21}^{(+)}l_{31}^{(+)})l_{11}^{(+)}, \qquad
  l_{43}^-=-l_{21}^{(-)}l_{11}^{(+)}l_{22}^{(+)}  , \\
    l_{ii}^{(+)} & = & l_{i^{ \prime} i^{ \prime}}^{(-)} ,
     \qquad  l_{ii}^{(-)}=(l_{ii}^{(+)})^{-1},  \qquad
     i=1, \ldots ,4.
  \end{eqnarray*}
Commutation relations in $ U_z(sp(2;j)) $ are given by
 \begin{eqnarray}
 [H_i,H_j ] = 0, & [H_i,X_j^{  \pm }]= \pm ( \alpha_i,
  \alpha_j ) X_j^{ \pm }, &  \nonumber \\
 & \left [X_i^+,X_j^- \right]= j^{2}_i [H_i]_q \delta_{ij} . &  \label{2}
 \end{eqnarray}
 \begin{equation}
(\alpha_i, \alpha_j )= -2 \varepsilon \left (\begin{array}{cc}
             1 & -1  \cr
            -1 &  2
              \end{array} \right ).
 \end{equation}
In the case $ j = \iota  $ we have obtained a new nonsemisimple algebra.


\begin{thebibliography}{99}
\bibitem{1}
N.A. Gromov, {\it Contractions and analytical continuations of classical
groups.
Unified approach}, { Syktyvkar, ëÏmi SC, 1990} (in Russian).

\bibitem{2}
N.Yu. Reshetikhin, L.A. Takhtajan, L.D. Faddeev, {\it  Algebra and analysis.}
(1989), Vol. 1, No.1, pp. 178-206.

\end{thebibliography}
\end{document}